\documentclass[12pt,a4paper]{article}

\title{A Phase Transition for the Metric Distortion
       of Percolation on the Hypercube}

\author{Omer Angel \and Itai Benjamini}
\date{June, 2003}

\usepackage{amsfonts,amsmath,amsthm}

\newtheorem{thm}{Theorem}
\newtheorem{lemma}[thm]{Lemma}
\newtheorem{prop}[thm]{Proposition}
\newtheorem{defn}[thm]{Definition}

\newcommand{\ep}{\epsilon}
\newcommand{\Z}{\mathbb Z}
\newcommand{\E}{\mathbb E}
\renewcommand{\P}{\mathbb P}

\begin{document}
\maketitle

\begin{abstract}
Let $H_n$ be the hypercube $\{0,1\}^n$, and let $H_{n,p}$ denote
the same graph with Bernoulli bond percolation with parameter
$p=n^{-\alpha}$. It is shown that at $\alpha=1/2$ there is a phase
transition for the metric distortion between $H_n$ and $H_{n,p}$.
For $\alpha<1/2$, asymptotically there is a map from $H_n$ to
$H_{n,p}$ with constant distortion (depending only on $\alpha$).
For $\alpha>1/2$ the distortion tends to infinity as a power of
$n$. We indicate the similarity to the existence of a non-uniqueness
phase in the context of infinite nonamenable graphs.
\end{abstract}

\section{Introduction}
%%%%%%%%%%%%%%%%%%%%%%

The hypercube $H_n$ is the graph with vertex set $\{0,1\}^n$ and
edges between two vectors differing in a single coordinate.
Bernoulli bond percolation on a graph $G$ with parameter $p$ is a
distribution on subgraphs $H\subset G$ where each edge is in $H$
with probability $p$ independently of all other edges. Site
percolation is defined similarly except that each vertex is in $H$
with probability $p$ and all edges with endpoints in $H$ are in $H$.
$H_{n,p}$ is the graph resulting from bond Bernoulli percolation
on $H_n$. We view a graph $G$ as a metric space on its vertex set
$V(G)$ with the induced graph metric, i.e.\ $d(a,b)$ is the length
of the shortest path from $a$ to $b$. Our interest is in the effect of
Bernoulli percolation on the geometry of the hypercube.

In \cite{hln} Hastad, Leighton and Newman show that if $p$ is
taken to be constant, then for large $n$ the distortion between
$H_n$ and $H_{n,p}$ is at most a constant with probability tending
to 1. If a graph is deemed to represent a computer network, then they were
interested in the latency of embedding $H_n$ in $H_{n,p}$, which measures
the loss of computational power resulting from eliminating a subset of the
processors (vertices) or communication channels (edges). The latency of an
embedding between graphs is slightly stronger then metric distortion.

In the classic paper \cite{aks} of Ajtai, Koml\'os and
Szemer\'edi, bond percolation on the hypercube is analyzed. It is
shown that if $p=an^{-1}$ with $a<1$, then all connected
components of $H_{n,p}$ have size polynomial in $n$ (hence
logarithmic in the size of $H_n$), while if $a>1$ there is a
single giant component of size linear in the size of $H_n$, and of
diameter of order $n$. When all connected components are small it
is evident that there is no map with a small distortion between
the cube and any of the components. However, when a giant
component exists it is a priori possible that the distortion
between the cube and the giant component is constant. Thus there is a gap
between $p=n^{-1}$ and $p$ constant where it is possible that the
distortion is small.

To make things more precise, we use the following measure of the distortion
between metric spaces.

\begin{defn}
The {\em distortion} of a map $f$ between metric spaces $(X,d_X),(Y,d_Y)$
is given by $D(f) = D_+(f)/D_-(f)$
\begin{gather*}\label{eq:dist}
D_+(f) = \max\left\{ 1, \sup_{a,b\in X} \frac{d_Y(f(a),f(b))}{d_X(a,b)}
             \right\} , \\
D_-(f) = \inf_{a,b\in X} \frac{\max\{1,d_Y(f(a),f(b))\}}{d_X(a,b)} .
\end{gather*}
The metric distortion between $X$ and $Y$, $D(X, Y)$, is
\[
\inf_{f:X\to Y} D(f) .
\]
\end{defn}

Note that since we will deal with functions that are not
necessarily 1 to 1, the definition accommodates the possibility that
two points $a,b$ are mapped to the same point in $Y$. This does
not constitute a significant change from other definitions. It should be
noted though that this definition is not symmetric, as the functions do not
need to have a full image.

Two graphs (or metric spaces) with a small metric distortion
between them  have several similar properties. For more
information on metric distortion see \cite{li} and reference
there. Note that the distance between disjoint connected
components of a graph is infinite, so only maps that preserve
connected components are of interest.

Note that (for bond percolation) a graph $G$ and $G_p$ have the same
vertex sets, so there is a natural map between the two. However,
there is no guarantee (and generally it is not the case) that this
map gives a minimal distortion or even gets close to doing that.
However, in the case considered here, it will turn out that when
the metric distortion is small, there is a map that approximates
the identity map which gives a small distortion.

\begin{thm}\label{thm:main}
Let $p=n^{-\alpha}$.
\begin{enumerate}
\item If $\alpha<1/2$ then there is some constant $c=c(\alpha)$ such that
\[
\P(D(H_n,H_{n,p}) < c) \to 1
\]
as $n\to\infty$.
\item If $\alpha>1/2$ then there is some constant $\beta=\beta(\alpha)>0$
such that
\[
\P(D(H_n,H_{n,p}) < n^\beta) \to 0
\]
as $n\to\infty$ .
\end{enumerate}
\end{thm}

A similar theorem holds for site percolation, as well as for mixed
percolation (with threshold at $p_bp_s=n^{-1}$). The essentially
identical proof will be omitted.

It follows from the definition of the latency of an embedding (see
\cite{hln}) that for $p=n^{-\alpha}$ with $\alpha<1/2$ the
latency is polynomial in $n$. This follows from
Theorem~\ref{thm:main}, as the latency is at most
$n^{D(f)}$. It is plausible (though not proved here) that when
$\alpha>1/2$ the latency is super-polynomial in $n$.

The phase transition in this context is linked below to the
$p_c<p_u$ conjecture in context of infinite nonamenable graphs. In
this light it is interesting to compare the behavior of the metric
distortion of percolation in the hypercube to the behavior in
other natural graphs. For the $n$-box in $\Z^d$ the distortion
between the box and the giant component is $O(\log n)$. While not
constant, this is small in relation to the diameter and there is
no additional phase transition beyond the percolation threshold.
For percolation on the complete graph, the metric distortion is
given by the diameter of the giant component. On the hypercube the
giant component has diameter $\theta(n)$ a.s.\ for any
$p>p=an^{-1}$ with $a>1$.

Regarding the proof, the main idea for $p=n^{-\alpha}$ with
$\alpha<1/2$, is that two vertices that are nearby in $H$ are
likely to stay nearby in $H_{n,p}$. Thus the identity map is almost
good enough, as it gives a constant distortion for most pairs of
vertices. To find a map that gives a constant distortion requires
a somewhat finer choice of the map, though the distance from $x$
to $f(x)$ is still bounded.

When $\alpha>1/2$ the proof is based on showing that while the
cube contains ``a lot" of geodesic cycles, the giant component for
$\alpha>1/2$ is locally very treelike and contains much fewer geodesic
cycles and thus distortion is generated.

In this regime there is a constant probability that two adjacent
vertices that are both in the giant component are not joined by a
short path in $H_{n,p}$, but will be at distance proportional to
the diameter. Of course, since for the cube (and much more
generally, \cite{aks} \cite{abs}) the giant component is known to
be unique, the distance between the two vertices is finite.
However, there is no longer a direct correlation between the
metric of the cube and that of the giant component. As giant
components are the finite graph analogue of infinite percolation
clusters, this behavior can be seen as the finite graph analogue
of non-uniqueness of infinite clusters. For graphs where $p_u$ ---
the threshold for uniqueness of the infinite cluster --- is
strictly greater than $p_c$ there is a regime where two nearby
vertices have a positive probability of being in disjoint infinite
components . That strict inequality was established in many cases
of nonamenable infinite graphs, see e.g. \cite{l} \cite{bs}. We
believe that the distortion in this case is indeed linear in the
diameter. However, since we have not eliminated possible maps that
are far from preserving the geometry we only give a lower bound of
$n^\beta$ for some $\beta$.

\medskip

Thus we propose that an analogue in finite graphs of having
non-unique infinite components is that the metric in the giant
component is very different from the original graph metric. In
particular, with probability bounded away from $0$, neighbors in
the graphs are at distance of the order of the diameter in the giant
component. For more on percolation on finite graphs see \cite{abs}.

An additional aspect of the phase transition at $\alpha=1/2$ is the routing
complexity. The routing problem is to find a path in $H_{n,p}$ between two
given vertices, preferably a short one. The routing problem may be defined
under one of several models, the strictest of which involves starting at
$x$ and only being allowed to query edges incident on vertices that have
already been reached from $x$. We refer to this as the local model.

For $p=n^{-\alpha}$ with $\alpha<1/2$, it turns out that there is a
polynomial time (in $n$) algorithm for routing in the local model,
outputting a path of length $O(n)$ whenever $x,y$ are in the same component.
This is a consequence of the fact that distance between neighboring
vertices is typically bounded, and so to get from $x$ to a nearby $y$ only
a ball of constant radius (and polynomial volume) needs to be explored.
On the other hand when $\alpha>1/2$ no polynomial algorithm exists.
For more details on complexity of routing in the hypercube and other
scenarios see \cite{abow}.

\medskip

We now proceed to present the proof of Theorem~\ref{thm:main}. The
proof is separated into the super-critical and the sub-critical
cases, each requiring a different approach.

\section{ The Super-Critical Case: $\alpha<1/2$  }
%%%%%%%%%%%%%%%%%%%%%%%%%%%%%%%%%

To illustrate the proof we first show a weaker result, namely that the
typical distortion of the identity map is constant.

\begin{prop}\label{prop:typical}
Let $p=n^{-\alpha}$ with $\alpha<1/2$. For some $l=l(\alpha)$, with
probability tending to 1 the percolation distance between two neighbors in
$H_n$ is T most $2l+1$.
\end{prop}

\begin{proof}[Sketch of Proof]
We consider paths of length $2l+1$ between them making $l$ steps in distinct
coordinates, making a step in the coordinate that differs $x$ from $y$, and
then retracing the first $l$ steps to reach $y$. The number of such paths
is roughly $n^l$ and each is open with probability $p^{2l+1}$. The
expectation of $X$: the number of open paths is
$n^{(1-2\alpha)l-\alpha}(1+o(1))$, and for some $l$ the exponent is
positive. 

A pair of paths is unlikely to have a large intersection. The second moment
of $X$ can be approximated by the square of the expectation. The largest
error term comes from pairs of paths that intersect only in their first
step, so that $\E X^2 =(\E X)^2 (1+ n^{2\alpha-1}+o(n^{2\alpha-1}))$. The
second moment method (see \cite{as}) yields that
\[
\P(X=0) \le 1-\frac{(\E X)^2}{\E X^2} = (1+o(1))n^{2\alpha-1} \to 0 .
\]
\end{proof}

Similarly, the idea in the following proof is that there are many paths of
some length between nearby vertices, and that one of them is very likely to
be open. However, in order to get an exponential bound on the probability
of failure we choose the set of paths we consider carefully.

Fix $l$ so that $(1-2\alpha)l>9\alpha$, which is possible as $\alpha<1/2$.
Set $m=\lfloor \frac{n-1}{l+2} \rfloor$ and fix throughout this section some
partition of the coordinates into disjoint sets of size $m$ denoted
$A,B,C_1,\ldots,C_l$ (with at least one coordinate left unused). A vertex
in $H_{n,p}$ 
will be called {\em good} if it has at least $2m$ vertices at distance 2
(in $H_{n,p}$) differing only in $A$-coordinates. We first show that good
vertices are dense: 

\begin{lemma}
Let $p=n^{-\alpha}$ for some $\alpha<1/2$. With probability
tending to 1, every vertex $v$ has is some good vertex $u$ in
$H_{n,p}$ differing from $v$ only in a single $B$-coordinate.
\end{lemma}

\begin{proof}
First we show that the probability that a vertex is good tends to 1 as
$n\to\infty$, with $\alpha$ fixed. To this end we count open paths of
length 2 using $A$-coordinates starting at $v$. If there are more that $4m$
such paths, then they connect $v$ to at least $2m$ vertices and $v$ is good.
The probability that $v$ is incident on at least
$m \cdot p/2 = c_1 n^{1-\alpha}$ open $A$-edges tends to 1 as $n\to\infty$.
On this event there are at least $(m-1) c_1 n^{1-\alpha}$ edges that extend
these to paths of length 2. With probability tending to 1 at least fraction
$p/2$ of these are open giving a total of at least $c_2 n^{2-2\alpha}$ open
paths. Since for $n$ large enough, $c_2 n^{2-2\alpha} > 4m$, the
probability that $v$ is good tends to 1.

Fix some $v$ and let $u_1,\ldots,u_m$ be $u$'s neighbors along $B$
coordinates. Since the $u_i$'s all differ in their $B$-coordinates, the
events of the $u_i$'s being good are independent. For large enough $n$ each
$u_i$ is good with probability at least $1-3^{-l}$ and then the probability
that none of the $u_i$'s is good is at most $3^{-ml} \ll 2^{-n}$. Thus, a
union bound over all vertices $v$ shows that with probability tending to 1 
every $v$ has a good neighboring $u$.
\end{proof}

\begin{proof}[Proof of Theorem~\ref{thm:main}(1)]
We construct a map $f$ that takes each vertex $x$ to some good
neighbor of $x$ by changing a single $B$ coordinate. By the
above Lemma, with probability tending to 1 such a map exists.
Since $d(x,f(x))=1$, it follows that $d(f(x),f(y)) \geq d(x,y)-2$
and so $D_-(f)>1/3$.
To show that distances do not increase by more then some constant
factor, it suffices to show that with probability tending to 1
the percolation distance between $f(x)$ and $f(y)$ for all
neighboring $x,y$ is bounded by $2l+9$.

The distance between $f(x)$ and $f(y)$ is at most 3. First, fix some
coordinate $e$ where $f(x)$ and $f(y)$ differ. In the remainder there are
some slight variations according as to which set $e$ is part of. Since
$f(x)$ is good it has $2m$ vertices at percolation distance $2$, differing 
only in the $A$-coordinates. Let $x_1,\ldots,x_m$ be $m$ of those that do
not differ from $f(x)$ in the $e$ coordinate (If $e\in A$ there are at
least $m$ such neighbors; Otherwise all $2m$ are suitable).
Similarly define $y_1,\ldots,y_n$. Note that $x_i$ differs from $y_i$ in
at most 7 coordinates.

We look for an open path from $x_i$ to $y_i$ of length $2l+9$ for
some $i$, and bound the probability that there is no such
path. The second moment method is used to bound the probability of
having no path for some $i$. The sets of paths we will consider
are chosen so that distinct $i$'s will be independent.

To each $i$ we associate a unique coordinate $b_i\in B$. If $e\in B$ we use
instead of $e$ the remaining coordinate not in any of the sets. For a
sequence $\mathbf c=(c_1,\ldots,c_l)$ of coordinates with $c_k\in C_k$ we
define the path $P_i(\mathbf j)$ from $x_i$ to $y_i$ as the path
making $l$ steps in the coordinates given by $\mathbf c$, making a
step in the coordinate $b_i$, making steps in the coordinates where
$x_i$ and $y_i$ differ, and then retracing the $l+1$ steps given by
$i$ and $\mathbf c$ in reversed order, thereby arriving at $y_i$.
If $e\in C_k$ for some $k$ then we again replace it by the unused
coordinate.

We now claim that paths $P_i(\mathbf c)$ and $P_{i'}(\mathbf c')$
for $i\neq i'$ are disjoint. Indeed, suppose that some $v$ is in some
$P_i(\mathbf c)$. According to the $e$ coordinate of $v$ we can say whether
$v$ is in the first or second half of the path. Assume w.log.\ it is the
first. By the distance from $f(x)$ we can determine the exact position of
$v$ in the path. Finally, either $v$ differs from $f(x)$ in the $b_i$
coordinate, and then we can say what $i$ is, or it differs from $f(x)$ in
exactly 2 $a$ coordinates and then we can say which is $x_i$. In any case,
we can determine $i$ from $v$.

\medskip

Let $X$ be the number of open paths from $x_i$ to $y_i$.
\[
\E X = p^{2l+9} m^l = cn^{(1-2\alpha)l-9\alpha} ,
\]
where $c$ is a constant depending only on $l$. For $l=l(\alpha)$ as
specified above the the exponent is positive, and so the expected number
of paths tends to infinity. The second moment of $X$ is given by
\begin{eqnarray*}
\E X^2
  & = & \sum_{P,P'} p^{4l+18-|P\cap P'|}        \\
  &\le& (\E X)^2 \left(\sum_k^{l-1} p^{-2k} m^{-k} + m^{-l}p^{-2l-9} \right)\\
  &\le& (\E X)^2 (1+c n^{2\alpha-1}) .
\end{eqnarray*}

And so the probability of having $X>0$ tends to 1
as $n\to\infty$. In particular, for large enough $n$, the
probability of having no open path for any of the $i$ is 
less then $(3^{-l})^m$. In this case there is an open path
of length $2l+13$ between $f(x)$ and $f(y)$. Since the number of
pairs $x,y$ is $n 2^n \ll 3^{lm}$, a union bound shows that asymptotically
almost surely distances increase by a factor of no more than
$2l+13$.
\end{proof}

\section{The Sub-Critical Case: $\alpha>1/2$}
%%%%%%%%%%%%%%%%%%%%%%%%%%%%%%%%%%%%%%%%%%%%%

For two vertices $x,y$ at some fixed distance from each other, the
number of paths of length $l$ between them grows like $n^{l/2}$,
and so the probability of having an open path tends to 0 as
$n\to\infty$. This shows that there is no map with bounded $d(x,f(x)$ with
small distortion. However, the image of a map need not cover all pairs
$x,y$, and there will be some pairs where paths will exist.

A {\em geodesic cycle} in a graph is an isometric embedding of a cycle
in the graph. The cube has a great number of geodesic cycles
passing through any given vertex. These can be found by making any
number of steps in distinct coordinates, and then repeating the
same sequence of coordinates to return to the starting point.

A geodesic cycle (as well as a geodesic path) is a structure that
is roughly preserved by a map with a small distortion. The key
idea of the proof is to show that the percolated cube does not
have many cycles, restricting the possible maps.

\begin{lemma}
If a map from $G$ to $H$ has distortion $D$, and $G$ has a geodesic cycle
$C$ of length $2l$ passing through a vertex $v$, then there is a simple
(self avoiding) cycle $C'\subset H$ of length $l' \in [l/2D,Dl]$ so that
$f(v)$ is at distance at most $2D^2$ from $C'$.
\end{lemma}

The key idea is to look at a cycle in $H$ that is the ``image'' of
$C$, and remove any small loops it has to extract a simple cycle.

\begin{proof}
Recall that the distortion is $D=D_+/D_-$.
Suppose the vertices of $C$ are $v_0,v_1,\ldots,v_l$ with
$v_0=v_l=v$, and their images in $H$ are $x_i=f(v_i)$. Since
$d(x_i,x_{i+1})\le D_+$, there are simple paths of length at most
$D_+$ between them, and these paths may be joined to form a (not
necessarily simple) cycle $C_0$ in $H$ of length at most $D_+l$,
passing through all the $x_i$'s in order. The paths from $x_i$ to
$x_{i+1}$ for various $i$'s will be denoted as arcs.

Next we extract from $C_0$ a simple cycle by a erasing loops it
includes. Note that the procedure we use is not the standard
procedure known as loop erasure. Suppose some vertex $u$ appears
twice in $C_0$. Since $u$ cannot be used twice in the same arc, it
appears in distinct arcs. If $u$ appears in the arc
$(x_i,x_{i+1})$ and in $(x_j,x_{j+1})$ then there is a path
of length at most $2D_+$ between $x_i$ and $x_j$. However,
$d(x_i,x_j) \ge D_- d(v_i,v_j)$, and so $d(v_i,v_j) \le 2D$. Since
the $v_i$'s form a geodesic cycle, breaking $C_0$ at the two
occurrences of $u$ results in one segment containing all but at
most $2D$ of the $x_i$'s.

Since the above holds for any 2 appearances of $u$ in $C_0$, for each
vertex $u\in C_0$ there is a single segment of $C_0$ starting and ending
at $u$, including all appearances of $u$, and including at most $2d$ of the
$x_i$'s. We now repeatedly remove from $C_0$ the longest such segment
corresponding to any repeated vertex. This choice guarantees that if a loop
starting and ending at $u$ is removed, then $u$ itself will not be removed
later. Since at each step we remove at most $2D$ of the $x_i$'s and create
a vertex that cannot be removed at a later stage, eventually a simple cycle
is left of length at least $l/2D$.

Finally, if $x_0=f(v)$ is not in the remaining cycle, it is in a loop that
was removed from the cycle at some stage. Since the length of each
such loop is at most $2D D_+$, it follows that $f(v)$ is close to
the cycle.
\end{proof}

\begin{lemma}\label{lem:no_loops}
Let $p=n^{-\alpha}$ with $\alpha>1/2$, and $\gamma<2\alpha-1$. For
large enough $n$, the probability that an vertex $v$ is included
in an open simple cycle in $H_{n,p}$ of length $2l \in [2n^\beta,
2n^{\gamma}]$ is at most $n^{1+n^\beta(\beta+1-2\alpha)}$.
Furthermore, the probability that $v$ is within distance $\delta$ of
such a cycle is at most $n^{\delta+1+n^\beta(\beta+1-2\alpha)}$.
\end{lemma}

\begin{proof}
First we estimate the number of simple cycles of length $2l$ in
$H_n$ that include $v$. Since a cycle makes an even number of
moves in each coordinate, the steps of a cycle of length $2l$ can
be partitioned into $l$ pairs of (non consecutive) moves in the
same coordinate. The number of partitions of $2l$ elements into
pairs is $(2l-1)!!$. The number of ways to choose a
coordinate for each of the pairs is $n^l$. The probability that a
cycle of length $2l$ is open is $p^{2l}=n^{-2\alpha l}$, so the
probability of $v$ being in an open cycle of length $2l$ is at
most $(2l-1)!!(n^{1-2\alpha})^l$.

This bound is decreasing in $l$ as long as
$2l<n^{2\alpha-1}$, so when summing over $l$, for large enough $n$ the
first term is the largest:
\[
\sum_{l=n^\beta}^{n^\gamma} (2l-1)!!(n^{1-2\alpha})^l
 < n (2n^\beta-1)!! (n^{1-2\alpha})^{n^\beta}
 < n^{1 + n^\beta(\beta+1-2\alpha)} .
\]
Finally, the number of vertices within distance $\delta$ of $v$ is at most
$n^\delta$, and the second claim follows from the first.
\end{proof}

\begin{proof}[Proof of Theorem~\ref{thm:main}(2)]
Fix positive constants $\beta,\gamma$ so that $\beta+\gamma<2\alpha-1$
and $\gamma>3\beta$. Assume that there is a map $f: H_n \to
H_{n,p}$ with distortion at most $n^\beta$. Since vertices at
distance greater than $n^\beta$ are mapped to distinct vertices,
the range of $f$ has size at least $2^n n^{-n^\beta}$. Since every vertex
in $H_n$ is in a geodesic cycle of length $n^\gamma$, every
vertex in the range is at distance at most $2n^{2\beta}$ from a
simple cycle of length $l \in [n^{\gamma-\beta},2n^{\gamma+\beta}]$.

However, Lemma~\ref{lem:no_loops} bounds the probability that a
vertex has this property, and thus the expected number of such
vertices. The probability of having at least $2^n n^{-n^\beta}$
vertices close to such cycles is at most
\[
\frac
{2^n n^{2n^{2\beta}+1+(\gamma-\beta+1-2\alpha)n^{\gamma-\beta}}}
{2^n n^{-n^\beta}}
 = n^{2n^{2\beta}+n^\beta +1+ (\gamma-\beta+1-2\alpha)n^{\gamma-\beta}}.
\]
This tends to 0 as $n\to \infty$ since the dominant term in the
exponent is $n^{\gamma-\beta}$ with a negative coefficient.

Thus the probability of having a map with distortion $n^\beta$
tends to 0. The above constraints allow taking
$\beta=(2\alpha-1)/4-\ep$ .
\end{proof}

\section{Problems}
%%%%%%%%%%%%%%%%%%

\begin{itemize}

\item
What is the critical window? What happens if $p=an^{-1/2}$? Is
there some $a_c$ so that above $a_c$ the distortion is constant or
perhaps logarithmic in $n$ while below it the distortion is
logarithmic in $n$ or even larger? What do the above proofs yield
in this case?

\item
When $\alpha>1/2$, what is the true behavior of the distortion? Since the
diameter of $H_n$ is $n$ the distortion can not be larger, but the proof
only gives a lower bound of $n^\beta$ for any $\beta<(2\alpha-1)/4$. 

\item
When $\alpha<1/2$, what is the true distortion? The proof of
Theorem~\ref{thm:main} gives a bound $D_+ = O((1/2-\alpha)^{-1})$. The
proof of Proposition~\ref{prop:typical} suggests a bound of
$D_+ = 1+\lfloor \frac{\alpha}{1/2-\alpha}\rfloor$. We believe this is
the correct uniform bound on the distortion.

\item
What is the minimal latency of an embedding of $H_n$ in $H_{n,p}$ when
$\alpha<1/2$? Show that it is super-polynomial in $n$. Does it grow
exponentially in $n$?

\item
Let $G$ be a vertex transitive graph. Perform $p$-bond percolation
on $G$. Assume that the  median distance between neighbors in $G$,
after percolation is $k$, show that the metric distortion between
$G$ and the percolation giant component is at least $k$.

\item
In particular, assume a vertex transitive graph $G$ has girth $g$,
is it true that the distortion between $G$ and a percolation
subgraph of $G$, with $p$ bounded away from $1$, is a.s.\ at least
$g$?

\end{itemize}

\bigskip

\noindent{\bf Acknowledgments.}
Thanks to Eran Ofek and Udi Wieder for useful discussions.

\end{document}